# The Balakrishnan Alpha Skew Laplace Distribution: Properties and Its Applications

Sricharan Shah, Partha J. Hazarika and Subrata Chakraborty

Department of Statistics, Dibrugarh University, Dibrugarh, Assam, 786004

charan.shah90@gmail.com, parthajhazarika@gmail.com, and subrata_stats@dibru.ac.in

**Abstract**

In this study by considering Balakrishnan mechanism a new form of alpha skew distribution is proposed and properties are investigated. The suitability of the proposed distribution has tested at the end with appropriate data fitting experiment and then comparing the values of Akaike Information Criterion (AIC), and Bayesian Information Criterion (BIC) with the values of some other related distributions. Likelihood ratio test is carried out for nested models, that is, for Laplace and the proposed distributions.

## 1. Introduction

The study of skew symmetric distribution was first introduced with the skew-normal distribution by Azzalini (1985) which is nothing but a natural extension of normal distribution derived by adding an additional parameter to introduce asymmetry. Its density function is given by

$$f_Z(z;\lambda) = 2g(z)G(\lambda z); \quad -\infty < z < \infty \tag{1}$$

where, $g(.)$ and $G(.)$ are the probability density function (pdf) and cumulative distribution function (cdf) of standard normal distribution respectively, and $\lambda \in R$, the skewness parameter. Denote it by $Z \sim \text{SN}(\lambda)$. Since then a lot of study have been carried out and a number of extensions and generalizations of this distribution were proposed and studied (for details see Chakraborty and Hazarika, 2011).

As a discussant in Arnold and Beaver (2002), Balakrishnan (2002) proposed the generalization of the skew normal density and studied its properties and the pdf of the same is

$$f_Z(z;\lambda,n) = g(z)[G(\lambda z)]^n / C_n(\lambda) \tag{2}$$

where, $n$ is positive integer and $C_n(\lambda) = E[G^n(\lambda U)]$, $U \sim N(0,1)$.

Huang and Chen (2007) developed the generalized skew-symmetric distributions by introducing the concept of skew function $G(.)$ instead of cdf in equation (1) where, $G(.)$ is a Lebesgue measurable function satisfying, $0 \leq G(z) \leq 1$ and $G(z) + G(-z) = 1$, $z \in R$, almost everywhere. Obviously, by choosing different skew functions, one can construct many numbers of skewed distributions (see Elal Olivero 2010; Harandi and Alamatsaz 2013; Hazarika and Chakraborty 2014; Chakraborty et al. 2012, 2014, and 2015 etc).

In this article *first*, a new class of skew distributions which is flexible enough to can be both unimodal and bimodal is proposed by considering Balakrishnan methodology of equation (2) for $n = 2$. *Second*, we apply this distribution in real life applications for better fitting by comparing this new distribution with the other known distributions.

The rest of this paper is organized as follows: In Section 2, we introduced a new form of alpha skew Laplace distribution and discuss its distributional properties. Section 3 discussed some extensions and generalization of this distribution. The parameter estimation of this distribution and its applications in real life is produced in Section 4. Finally, conclusions are provided in last section.

## 2. A New Alpha Skew Laplace Distribution

**Definition 1:** A random variable $Z$ is distributed as a Balakrishnan alpha skew Laplace distribution, denoted by $BASLa_2(\alpha)$, if its pdf is of the form

$$f_Z(z;\alpha) = \frac{[(1-\alpha z)^2 + 1]^2}{C_2(\alpha)} \frac{1}{2} e^{-|z|} \qquad (3)$$

Where, $C_2(\alpha) = 4(1 + 4\alpha^2 + 6\alpha^4)$ and $\alpha$ is the skewness parameter.

Furthermore, when $\alpha = 0$, we get the standard Laplace distribution.

i) When $\alpha \to \pm\infty$, we get the Bimodal Laplace distribution (see Hazarika and Chakraborty, 2014) given by

$$f(z) = \frac{z^4}{48} e^{-|z|}$$

and is denoted by $Z \sim BLa(4)$.

ii) If $Z \sim BASLa_2(\alpha)$ then $-Z \sim BASLa_2(-\alpha)$.

The plots of the pdf of $BASLa_2(\alpha)$ for different choices of the parameter $\alpha$ are plotted in figure 1.

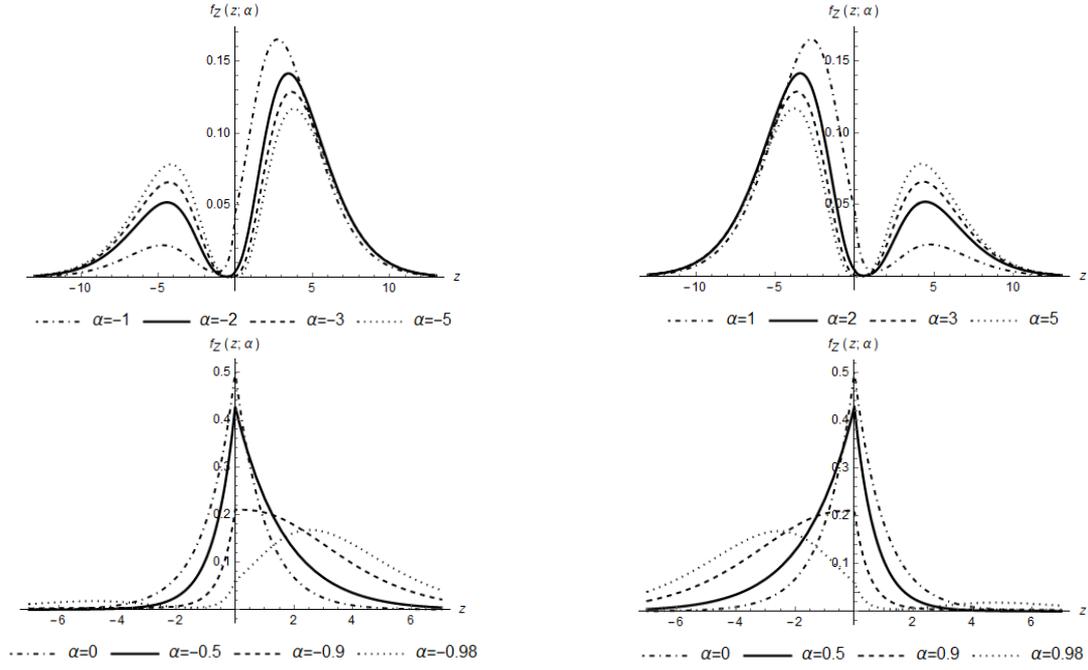

**Figure 1:** Plots of pdf of $BASLa_2(\alpha)$

**Proposition 1:** Let $Z \sim BASLa_2(\alpha)$. Then the cdf $Z$ of is given by

$$F_Z(z;\alpha) = \begin{cases} \dfrac{e^z\left(\begin{array}{l}4+\alpha(8-8z+8\alpha(2-2z+z^2))-4\alpha^2(-6+6z-3z^2+z^3)+\\ \alpha^3(24-24z+12z^2-4z^3+z^4))\end{array}\right)}{2C_2(\alpha)}; & z \leq 0 \\[2em] \dfrac{e^{-z}\left(\begin{array}{l}-4+8\alpha(z+1)-8\alpha^2(2+2z+z^2)+4\alpha^3(6+6z+3z^2+z^3)-\\ \alpha^4(24+24z+12z^2+4z^3+z^4)+2e^zC_2(\alpha)\end{array}\right)}{2C_2(\alpha)}; & z > 0 \end{cases}$$

(4)

**Proof:** $F_Z(z;\alpha) = \dfrac{1}{2C_2(\alpha)} \int_{-\infty}^{z} [(1-\alpha z)^2 + 1]^2 \, e^{-|z|} \, dz$

$= \dfrac{1}{2C_2(\alpha)} \left[ \int_{-\infty}^{z} 4e^{-|z|} \, dz - 8\alpha \int_{-\infty}^{z} z e^{-|z|} \, dz + 8\alpha^2 \int_{-\infty}^{z} z^2 e^{-|z|} \, dz - 4\alpha^3 \int_{-\infty}^{z} z^3 e^{-|z|} \, dz \right.$

$\left. + \alpha^4 \int_{-\infty}^{z} z^4 e^{-|z|} \, dz \right]$

**Case1:** When $z \leq 0$, then we get

$F_Z(z;\alpha) = \dfrac{1}{2C_2(\alpha)} \left[ 4(e^z) - 8\alpha(e^z(z-1)) + 8\alpha^2(e^z(2-2z+z^2)) - \right.$

$\left. 4\alpha^3(e^z(-6+6z-3z^2+z^3)) + \alpha^4(e^z(24-24z+12z^2-4z^3+z^4)) \right]$

On simplification we get the result in equation (5) for $z \leq 0$.

**Case2:** When $z > 0$, then we get

$$F_Z(z;\alpha) = \frac{1}{2C_2(\alpha)}\left[4\left(e^{-z}(-1+2e^z)\right) - 8\alpha\left(-e^{-z}(1+z)\right) + 8\alpha^2\left(e^{-z}(-2+4e^z-2z-z^2)\right) - \right.$$

$$\left. 4\alpha^3\left(-e^{-z}(6+6z+3z^2+z^3)\right) + \alpha^4\left(e^{-z}(-24+48e^z-24z-12z^2-4z^3-z^4)\right)\right]$$

On simplification we get the result in equation (5) for $z > 0$.

**Corollary 1:** In taking the limit $\alpha \to \pm\infty$ in equation (4), we get the cdf of $BLa(4)$ distribution as

$$F(z) = \frac{1}{48}\begin{cases} e^z(24-24z+12z^2-4z^3+z^4), & z < 0 \\ e^{-z}(-24+48e^z-24z-12z^2-4z^3-z^4), & z \geq 0 \end{cases}$$

The cdf is plotted in figure 2 for studying variation in its shape with respect to the parameter $\alpha$.

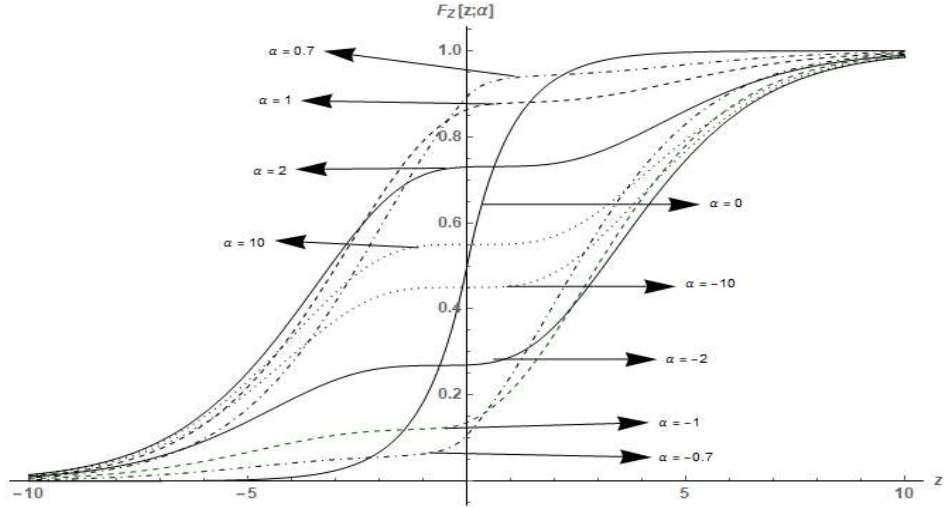

**Figure 2:** Plots of cdf of $BASLa_2(\alpha)$

**Proposition 2:** Let $Z \sim BASLa_2(\alpha)$. Then the mgf $Z$ of is given by

$$M_Z(t) = \frac{-(t^2-1)^4 - 4t\alpha(t^2-1)^3 - 4\alpha^2(t^2-1)^2(3t^2+1) - 24t\alpha^3(t^4-1) - 6\alpha^4(1+5t^4+10t^2)}{(t^2-1)^5(1+4\alpha^2+6\alpha^4)}$$

where, $-1 < t < 1$ (5)

*Proof:* $M_Z(t) = \dfrac{1}{2C_2(\alpha)}\int\limits_{-\infty}^{\infty} e^{tz}[(1-\alpha z)^2+1]^2 e^{-|z|}\,dz$

$= \dfrac{1}{2C_2(\alpha)}\left[\int\limits_{-\infty}^{\infty} 4e^{tz}e^{-|z|}\,dz - 8\alpha\int\limits_{-\infty}^{\infty} ze^{tz}e^{-|z|}\,dz + 8\alpha^2\int\limits_{-\infty}^{\infty} z^2 e^{tz}e^{-|z|}\,dz - 4\alpha^3\int\limits_{-\infty}^{\infty} z^3 e^{tz}e^{-|z|}\,dz\right.$

$\left. + \alpha^4\int\limits_{-\infty}^{\infty} z^4 e^{tz}e^{-|z|}\,dz\right]$

$$= \frac{1}{2\,C_2(\alpha)} \left[ 4\left( -\frac{2}{(t-1)(t+1)} \right) - 8\alpha\left( \frac{4t}{(t-1)^2(t+1)^2} \right) + 8\alpha^2\left( -\frac{4(1+3t^2)}{(t-1)^3(t+1)^3} \right) - \right.$$

$$\left. 4\alpha^3\left( \frac{48(t+t^3)}{(t-1)^4(t+1)^4} \right) + \alpha^4\left( -\frac{48(1+10t^2+5t^4)}{(t-1)^5(t+1)^5} \right) \right]$$

On simplification we get the desired result.

**Corollary 2:** In taking the limit $\alpha \to \pm\infty$ in equation (5), we get the mgf of $BLa(4)$ distribution as for $-1 < t < 1$

$$M_Z(t) = \frac{1+10t^2+5t^4}{(1-t)^5(t+1)^5}$$

**Proposition 3:** Let $Z \sim BASLa_2(\alpha)$. Then the $r^{th}$ moment of $Z$ is given by

$$E(Z^r) = \frac{\begin{pmatrix} 4(1+(-1)^r)\Gamma(1+r) + \alpha(8(-1+(-1)^r) + (2+r)\alpha(8(1+(-1)^r) + (3+r)\alpha(4(-1+ \\ (-1)^r) + (1+(-1)^r)(4+r)\alpha)))\Gamma(2+r) \end{pmatrix}}{2\,C_2(\alpha)}$$

(7)

**Proof:** $E(Z^k) = \dfrac{1}{2\,C_2(\alpha)} \int_{-\infty}^{\infty} z^r \left[(1-\alpha z)^2 + 1\right]^2 e^{-|z|} dz$

$$= \frac{1}{C_2(\alpha)} \left[ \int_{-\infty}^{\infty} 4z^r e^{-|z|} dz - 8\alpha \int_{-\infty}^{\infty} z^{r+1} e^{-|z|} dz + 8\alpha^2 \int_{-\infty}^{\infty} z^{r+2} e^{-|z|} dz - 4\alpha^3 \int_{-\infty}^{\infty} z^{r+3} e^{-|z|} dz \right.$$

$$\left. + \alpha^4 \int_{-\infty}^{\infty} z^{r+4} e^{-|z|} dz \right]$$

$$= \frac{1}{2\,C_2(\alpha)} \left[ 4\left\{(1+(-1)^r)\Gamma(1+r)\right\} - 8\alpha\left\{-(-1+(-1)^r)\Gamma(2+r)\right\} + 8\alpha^2\left\{(1+(-1)^r)\Gamma(3+r)\right\} - \right.$$

$$\left. 4\alpha^3\left\{-(-1+(-1)^r)\Gamma(4+r)\right\} + \alpha^4\left\{(1+(-1)^r)\Gamma(5+r)\right\} \right]$$

On simplification we get the result in equation (7).

Further, if we put $r = 1, 2, 3, 4$ in equation (7), we get

$$E(Z) = -\frac{4(\alpha+6\alpha^3)}{1+4\alpha^2+6\alpha^4}, \quad E(Z^2) = \frac{2(1+24\alpha^2+90\alpha^4)}{1+4\alpha^2+6\alpha^4}, \quad Var(Z) = \frac{2(1+20\alpha^2+96\alpha^4+216\alpha^6+540\alpha^8)}{(1+4\alpha^2+6\alpha^4)^2},$$

$$E(Z^3) = -\frac{48(\alpha+15\alpha^3)}{1+4\alpha^2+6\alpha^4}, \text{ and } E(Z^4) = \frac{24(1+60\alpha^2+420\alpha^4)}{1+4\alpha^2+6\alpha^4}.$$

**Remark 2:** The bounds for mean and variance can be obtained respectively by numerical optimization of $E(Z)$ and $Var(Z)$ with respect to $\alpha$ as $-2.58345 \le E(Z) \le 2.58345$ and $2 \le Var(Z) \le 30$, which can also be verified graphically in figure (3) and figure (4).

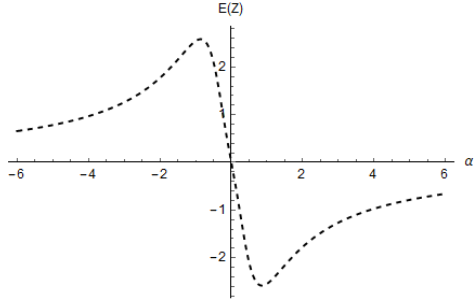
**Figure 3:** Plot of mean

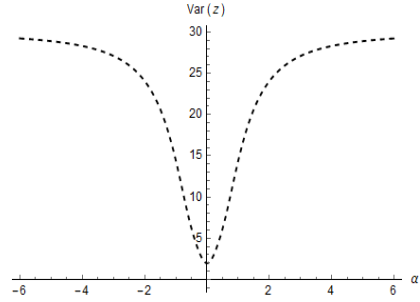
**Figure 4:** Plot of variance

And taking the limit $\alpha \to \pm\infty$ in the moments of $BASLa_2(\alpha)$ distribution, we get the moments of $BLa(4)$ distribution as

$$E(Z) \to 0, \quad Var(Z) \to 30$$

**Remark 4:** The skewness and kurtosis of $BASLa_2(\alpha)$ distribution are given by

$$\beta_1 = \frac{8\alpha^2(3+52\alpha^2+96\alpha^4-432\alpha^6-2700\alpha^8-6480\alpha^{10})^2}{(1+20\alpha^2+96\alpha^4+216\alpha^6+540\alpha^8)^3}$$

$$\beta_2 = \frac{6(1+48\alpha^2+566\alpha^4+3800\alpha^6+20820\alpha^8+77856\alpha^{10}+185544\alpha^{12}+246240\alpha^{14}+90720\alpha^{16})}{(1+20\alpha^2+96\alpha^4+216\alpha^6+540\alpha^8)^2}$$

**Remark 5:** The bounds for mean and variance can be obtained respectively by numerical optimization of $\beta_1$ and $\beta_2$ with respect to $\alpha$ as $0 \leq \beta_1 \leq 1.14182$ and $1.86667 \leq \beta_2 \leq 6.49587$, which can also be verified graphically in figure (5) and figure (6).

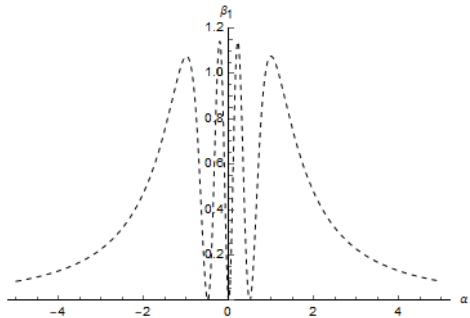
**Figure 5:** Plot of skewness

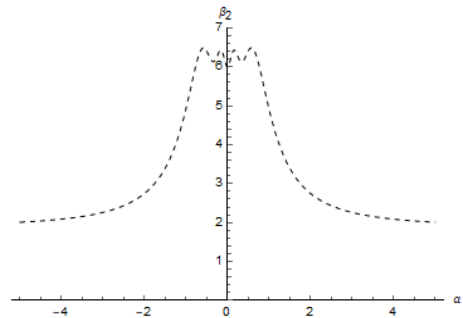
**Figure 6:** Plot of kurtosis

And, taking the limit $\alpha \to \pm\infty$ in the skewness and kurtosis of $BASLa_2(\alpha)$ distribution, we get the skewness and kurtosis of $BLa(4)$ distribution as

$$\beta_1 \to 0, \quad \beta_2 \to 1.8667$$

**Remark 6:** The pdf (3) can be written as

$$f_Z(z;\alpha) = \frac{4+8\alpha^2 z^2+\alpha^4 z^4}{C_2(\alpha)}\frac{1}{2}e^{-|z|} - \frac{8\alpha z+4\alpha^3 z^3}{C_2(\alpha)}\frac{1}{2}e^{-|z|} \qquad (8)$$

in which the 1st part is symmetric and the second part is asymmetric probability density function. We shall call a random variable $Z$ with the symmetric pdf i.e.,

$$f_1(z;\alpha) = \frac{4+8\alpha^2 z^2 + \alpha^4 z^4}{C_2(\alpha)} \frac{1}{2} e^{-|z|}; \quad \alpha \in R \qquad (9)$$

as the symmetric-component random variable of the model $BASLa_2(\alpha)$ and denote it by $SCBASLa_2(\alpha)$. Again, we note below that the pdf (9) can be expressed as a mixture between a Laplace density and a bimodal-Laplace density as

$$f_1(z;\alpha) = \frac{4}{C_2(\alpha)} \frac{1}{2} e^{-|z|} + \frac{8\alpha^2 z^2 + \alpha^4 z^4}{C_2(\alpha)} \frac{1}{2} e^{-|z|} \qquad (10)$$

Furthermore, we can generate random numbers for $BASLa_2(\alpha)$ by using the acceptance-rejection method in the following algorithm:

**Algorithm:**

Let $h(z)$ and $h_1(z)$ be the density function of $Z \sim BASLa_2(\alpha)$ and $H \sim SCBASLa_2(\alpha)$ random variables, respectively, and

$$T = \underset{z}{Sup} \frac{h(z)}{h_1(z)} = \frac{3+2\sqrt{2}}{3}$$

Now, we shall carry the following steps to simulate random numbers from $BASLa_2(\alpha)$ distribution as

(i) Generate a random variables $Y \sim SCBASLG_2(\alpha)$ and $U \sim \text{Uniform}(0,1)$, and they are independent.

(ii) If $U < \frac{1}{T} \frac{h(z)}{h_1(z)}$, and set $Z = H$; otherwise go back to step(i) and repeat the process.

## 3. Some Extensions and Generalization of $BASLa_2(\alpha)$ Distribution

### 3.1. A Two- parameter $BASLa_2(\alpha)$ Distribution

In this section, similar to the idea of Bahrami et al. (2009), we present the definition and some simple properties of two- parameter Balakrishnan alpha skew normal distribution.

**Definition 5:** A random variable $Z$ has a two-parameter Balakrishnan alpha skew normal distribution with parameters $\alpha_1, \alpha_2 \in R$, denoted by $TPBASN_2(\alpha_1, \alpha_2)$, if its pdf is

$$f(z;\alpha_1,\alpha_2) = \frac{\Phi^2(\alpha_1 z)\Phi^2(\alpha_2 z)}{C(\alpha_1,\alpha_2)} \frac{1}{2} e^{-|z|}; \quad z \in R \qquad (11)$$

Where,

$C(\alpha_1,\alpha_2) = 16[1+4\alpha_2^2 + 6\alpha_2^4 + 8\alpha_1(\alpha_2 + 6\alpha_2^3) + 48\alpha_1^3(\alpha_2 + 15\alpha_2^3) + 4\alpha_1^2(1+24\alpha_2^2 + 90\alpha_2^4) + 6\alpha_1^4\{1+60(\alpha_2^2 + 7\alpha_2^4)\}]$

$\Phi^2(\alpha_1 z) = [(1-\alpha_1 z)^2 + 1]^2$; and $\Phi^2(\alpha_2 z) = [(1-\alpha_2 z)^2 + 1]^2$

**Properties of $TPBASN_2(\alpha_1,\alpha_2)$:**

- If $\alpha_1 = \alpha_2 = 0$, then we get the standard Laplace distribution $f(z) = \frac{1}{2} e^{-|z|}$.

- If $\alpha_2 = 0$, then $Z \sim BASLa_2(\alpha_1)$ and if $\alpha_1 = 0$, then $Z \sim BASLa_2(\alpha_2)$.

- If $\alpha_1 = \alpha_2 = \alpha$, then

$$f(z;\alpha) = \frac{((1-\alpha z)^2 + 1)^4}{16(1 + 16\alpha^2 + 204\alpha^4 + 1440\alpha^6 + 2520\alpha^8)} \frac{1}{2} e^{-|z|} = BASLa_4(\alpha).$$

- If $\alpha_1 = \alpha_2 \to \pm\infty$, then $f(z) = \frac{z^8}{40320} \frac{1}{2} e^{-|z|}$.

- If $Z \sim TPBASN_2(\alpha_1, \alpha_2)$, then $-Z \sim TPBASN_2(-\alpha_1, -\alpha_2)$.

### 3.2. Balakrishnan Alpha-Beta Skew Laplace Distribution

In this section, using the work of Shafiei et al. (2016), we present the definition and some simple properties of Balakrishnan alpha-beta skew Laplace distribution.

**Definition 6:** If the density of the random variable $Z$ has pdf given by

$$f(z;\alpha,\beta) = \frac{[(1-\alpha z - \beta z^3)^2 + 1]^2}{C_2(\alpha,\beta)} \frac{1}{2} e^{-|z|} \quad ; \ z \in R \qquad (12)$$

Then, we say that $Z$ is distributed according to the Balakrishnan alpha-beta skew Laplace distribution with parameter $\alpha, \beta \in R$. where,

$C_2(\alpha,\beta) = 4(1 + 4\alpha^2 + 6\alpha^4 + 96\alpha\beta + 720\alpha^3\beta + 1440\beta^2 + 60480\alpha^2\beta^2 + 3628800\alpha\beta^3 + 119750400\beta^4)$

.We denote it by $Z \sim BAbSLa_2(\alpha,\beta)$.

**Properties of** $BAbSLa_2(\alpha,\beta)$ **:**

- If $\beta = 0$, then we get Balakrishnan alpha skew Laplace distribution and is given by

$$f(z;\alpha) = \frac{[(1-\alpha z)^2 + 1]^2}{4(1 + 4\alpha^2 + 6\alpha^4)} \frac{1}{2} e^{-|z|} .$$

- If $\alpha = 0$, then we get $f(z;\beta) = \frac{[(1-\beta z^3)^2 + 1]^2}{4(1 + 1440\beta^2 + 119750400\beta^4)} \frac{1}{2} e^{-|z|}$.

  This equation is known as Balakrishnan beta skew Laplace ($BbSLa_2(\beta)$) distribution.

- If $\alpha = \beta = 0$, then we get the standard Laplace distribution and is given by

$$f(z) = \frac{1}{2} e^{-|z|} .$$

- If $\alpha \to \pm\infty$, then we get the Bimodal Laplace ($BLa(4)$) distribution given by

$$f(z) = \frac{z^4}{48} e^{-|z|} .$$

- If $\beta \to \pm\infty$, then we get $f(z) = \dfrac{z^6}{720} e^{-|z|}$.

- If $Z \sim BAbSLa_2(\alpha, \beta)$, then $-Z \sim BAbSLa_2(-\alpha, -\beta)$.

### 3.3. Generalization of $BASLa_2(\alpha)$ Distribution

In this section, using the approach of Sharafi et al. (2017), we present the definition and some simple properties of generalization of alpha skew Laplace distribution.

**Definition 6:** If the random variable $Z$ has the pdf given by

$$f(z;\alpha,\lambda) = \frac{[(1-\alpha z)^2 + 1]^2}{C_2(\alpha,\lambda)} \frac{1}{2} e^{-|z|} \Psi(\lambda z); \quad z \in R \qquad (13)$$

Where, $C_2(\alpha,\lambda) = \dfrac{2}{\lambda}\left( (1+4\alpha^2 + 6\alpha^4)\lambda - \dfrac{2\alpha(2(1+6\alpha^2)\lambda^2 + (5+18\alpha^2)\lambda^3 + (1+3\alpha^2)\lambda^4(4+\sqrt{\lambda^2}))}{(1+\sqrt{\lambda^2})^4} \right)$

and $\Psi(\lambda z)$ is the cdf of Laplace distribution.

Then, we say that $Z$ is distributed according to the Generalized Balakrishnan alpha skew Laplace distribution with parameter $\alpha, \lambda \in R$. We denote it by $Z \sim GBASLa_2(\alpha, \lambda)$.

**Properties of $GBASLa_2(\alpha, \lambda)$:**

- If $\alpha = 0$, then we get $f(z;\alpha,\lambda) = e^{-|z|} \Psi(\lambda z)$.

- If $\lambda = 0$, then the function becomes Indeterminate.

- If $\alpha = \lambda = 0$, then $f(z) = \dfrac{1}{2} e^{-|z|}$.

- If $\alpha \to \pm\infty$, then $f(z;\alpha,\lambda) = \dfrac{z^4}{24} e^{-|z|} \Psi(\lambda z)$.

- If $Z \sim GBASLa_2(\alpha, \lambda)$, then $-Z \sim GBASLa_2(-\alpha, -\lambda)$.

### 3.4. The Log-Balakrishnan Alpha Skew Normal Distribution

In this section similar to the approach of Venegas et al. (2016), we present the definition and some simple properties of log-Balakrishnan alpha skew normal distribution.

Let $Z = e^Y$, then $Y = Log(Z)$, therefore, the pdf of Z is defined as follows:

**Definition 6:** If the random variable $Z$ has pdf given by

$$f(z;\alpha) = \frac{[(1-\alpha y)^2 + 1]^2}{C_2(\alpha)} \frac{1}{2z} e^{-|y|}; \quad z > 0 \qquad (14)$$

Then, we say that $Z$ is distributed according to the log-Balakrishnan alpha skew normal distribution with parameter $\alpha \in R$. Where, $y = Log(z)$ and $\varphi(y)$ is the pdf of the standard log-normal distribution. We denote it by $Z \sim LBASN_2(\alpha)$.

**Properties of** $LBASN_2(\alpha)$:

- If $\alpha = 0$, then we get $f(z) = \dfrac{1}{2z} e^{-|y|}$.

- If $\alpha \to \pm\infty$, then we get the Log-Bimodal Laplace $LBLa(4)$ distribution given by

$$f(z) = \frac{y^4}{24} \frac{1}{2z} e^{-|y|}.$$

- If $Z \sim LBASN_2(\alpha)$, then $-Z \sim LBASN_2(-\alpha)$.

## 4. Parameter Estimation and Applications

### 4.1. Location Scale Extension

By considering the location parameter ($\mu \in R$) and scale parameter ($\beta > 0$) extension, the pdf of $BASLa_2(\alpha)$ distribution can be written as

$$f_Z(z;\mu,\beta,\alpha) = \frac{\left[\left(1-\alpha\left(\dfrac{z-\mu}{\beta}\right)\right)^2 + 1\right]^2}{C_2(\alpha)} \frac{e^{-\frac{|z-\mu|}{\beta}}}{2\beta}; \quad z,\mu,\beta,\alpha \in R \qquad (15)$$

We denote it by $Z \sim BASLa_2(\alpha,\mu,\beta)$.

### 4.2. Maximum Likelihood Estimation

Consider a random sample $z_1, z_2, ..., z_n$ of size $n$ from the distribution of the random variable $Z \sim BASN_2(\alpha,\mu,\beta)$, so that the log-likelihood function for $\theta = (\alpha,\mu,\beta)$ is given by

$$l(\theta) = 2\sum_{i=1}^{n} \log\left[\left\{1-\alpha\left(\dfrac{z_i-\mu}{\beta}\right)\right\}^2 + 1\right] - n\log\beta - n\log 4 - n\log(1+4\alpha^2+6\alpha^4) + n\log\left(\dfrac{1}{2}\right) - \sum_{i=1}^{n}\dfrac{|z_i-\mu|}{\beta}$$

(16)

Using GenSA package in R, the MLEs of the parameters are obtained by numerically maximizing the log-likelihood function.

### 4.3. Real Life Applications

For the applicability of the proposed model, two real life datasets were analysed. The first dataset is related to N latitude degrees in 69 samples from world lakes, which appear in Column 5 of the Diversity data set in website: http://users.stat.umn.edu/sandy/courses/8061/datasets/lakes.lsp. The second dataset is related to the exchange rate data of the United Kingdom Pound to the United States Dollar from 1800 to 2003, which is available in the website http://www.globalfindata.com.

Further, we compared the proposed distribution $BASLa_2(\alpha,\mu,\beta)$ with the normal distribution $N(\mu,\sigma^2)$, the logistic distribution $LG(\mu,\beta)$, the Laplace distribution $La(\mu,\beta)$, the skew-normal distribution $SN(\lambda,\mu,\sigma)$ of Azzalini (1985), the skew-logistic distribution $SLG(\lambda,\mu,\beta)$ of Wahed and Ali (2001), the skew-Laplace distribution $SLa(\lambda,\mu,\beta)$ of Nekoukhou and Alamatsaz (2012), the alpha-skew-normal distribution $ASN(\alpha,\mu,\sigma)$ of Elal-Olivero (2010), the alpha-skew-Laplace distribution $ASLa(\alpha,\mu,\beta)$ of Harandi and Alamatsaz (2013), the alpha-skew-logistic distribution $ASLG(\alpha,\mu,\beta)$ of Hazarika and Chakraborty (2014). AIC and BIC are used for model selection.

**Table 1:** MLE's, log-likelihood, AIC and BIC for N latitude degrees in 69 samples from world lakes.

| Parameters → Distribution ↓ | $\mu$ | $\sigma$ | $\lambda$ | $\alpha$ | $\beta$ | log $L$ | AIC | BIC |
|---|---|---|---|---|---|---|---|---|
| $N(\mu,\sigma^2)$ | 45.165 | 9.549 | -- | -- | -- | -253.599 | 511.198 | 515.666 |
| $LG(\mu,\beta)$ | 43.639 | -- | -- | -- | 4.493 | -246.645 | 497.29 | 501.758 |
| $SN(\lambda,\mu,\sigma)$ | 35.344 | 13.69 | 3.687 | -- | -- | -243.036 | 492.072 | 498.774 |
| $SLG(\lambda,\mu,\beta)$ | 36.787 | -- | 2.828 | -- | 6.417 | -239.053 | 490.808 | 490.808 |
| $La(\mu,\beta)$ | 43 | -- | -- | -- | 5.895 | -239.248 | 482.496 | 486.964 |
| $ASLG(\alpha,\mu,\beta)$ | 49.087 | -- | -- | 0.861 | 3.449 | -237.351 | 480.702 | 487.404 |
| $SLa(\lambda,\mu,\beta)$ | 42.3 | -- | 0.255 | -- | 5.943 | -236.90 | 479.799 | 486.501 |
| $ASLa(\alpha,\mu,\beta)$ | 42.3 | -- | -- | -0.22 | 5.439 | -236.079 | 478.159 | 484.861 |
| $ASN(\alpha,\mu,\sigma)$ | 52.147 | 7.714 | -- | 2.042 | -- | -235.37 | 476.739 | 483.441 |
| $BASLa_2(\alpha,\mu,\beta)$ | 52.0 | -- | -- | 1.046 | 2.641 | **-232.138** | **470.277** | **476.979** |

**Table 2:** MLE's, log-likelihood, AIC and BIC for the exchange rate data of the United Kingdom Pound to the United States Dollar from 1800 to 2003.

| Parameters → Distribution ↓ | $\mu$ | $\sigma$ | $\lambda$ | $\alpha$ | $\beta$ | log $L$ | AIC | BIC |
|---|---|---|---|---|---|---|---|---|
| $N(\mu,\sigma^2)$ | 4.117 | 1.381 | -- | -- | -- | -355.265 | 714.529 | 721.165 |
| $SN(\lambda,\mu,\sigma)$ | 3.589 | 1.478 | 0.501 | -- | -- | -355.217 | 716.434 | 726.388 |
| $LG(\mu,\beta)$ | 4.251 | -- | -- | -- | 0.753 | -351.192 | 706.385 | 713.021 |
| $SLG(\lambda,\mu,\beta)$ | 5.36 | -- | -2.371 | -- | 1.018 | -341.391 | 688.782 | 698.736 |
| $La(\mu,\beta)$ | 4.754 | -- | -- | -- | 0.971 | -339.315 | 682.63 | 689.265 |
| $ASN(\alpha,\mu,\sigma)$ | 3.656 | 0.883 | -- | -3.504 | -- | -317.946 | 641.892 | 651.847 |
| $SLa(\lambda,\mu,\beta)$ | 4.861 | -- | 1.183 | -- | 0.439 | -302.564 | 611.128 | 621.082 |
| $ASLG(\alpha,\mu,\beta)$ | 3.764 | -- | -- | -2.025 | 0.403 | -301.963 | 609.927 | 619.881 |
| $ASLa(\alpha,\mu,\beta)$ | 4.861 | -- | -- | 0.539 | 0.677 | -301.443 | 608.885 | 618.840 |
| $BASLa_2(\alpha,\mu,\beta)$ | 4.871 | -- | -- | 0.296 | 0.603 | **-296.849** | **599.698** | **609.652** |

### 4.4. Likelihood Ratio Test

The likelihood ratio (LR) test is used to discriminate between $La(\mu,\beta)$ and $BASLa_2(\alpha,\mu,\beta)$ as they are nested models. For LR test, let us setup the Null hypothesis $H_0: \alpha = 0$, that is the sample is drawn from $La(\mu,\beta)$; against the alternative $H_1: \alpha \neq 0$, that is the sample is drawn from $BASLa_2(\alpha,\mu,\beta)$. For the Dataset 1, the value of LR test statistic (i.e., 14.22) which exceed the 99% critical value (i.e., 6.635). Thus, it supports the alternative hypothesis that the sampled data comes from $BASLa_2(\alpha,\mu,\beta)$, not from $La(\mu,\beta)$.

Again, for the Dataset 2, the value of LR test statistic (i.e., 84.932) which exceed the 99% critical value (i.e., 6.635). Thus, it supports the alternative hypothesis that the sampled data comes from $BASLa_2(\alpha,\mu,\beta)$, not from $La(\mu,\beta)$.

### 5. Conclusion

In this article a new alpha-skew-Laplace distribution which includes both unimodal as well as bimodal behavior is introduced and some of its basic properties are investigated. From the computation, it is examined that the proposed $BASLa_2(\alpha,\mu,\beta)$ distribution provides better fitting to the dataset under consideration in terms of all the criteria, namely log-likelihood, AIC and BIC. The plots of observed and expected densities presented above, also confirms our findings.

**References:**

1. Arnold, B. C., Beaver, R. J. (2002). Skewed multivariate models related to hidden truncation and/or selective reporting. *Test*, 11(1), 7-54.

2. Azzalini, A. (1985). A class of distributions which includes the normal ones. *Scandinavian Journal of Statistics*, 171-178.

3. Bahrami, W., Agahi, H., and Rangin, H. (2009). A two-parameter Balakrishnan skew-normal distribution. *J. Statist. Res. Iran* 6, 231-242

4. Balakrishnan, N. (2002). Discussion on "Skew multivariate models related to hidden truncation and/or selective reporting" by B.C.Arnold and R.J.Beaver. *Test* 11:37-39

5. Chakraborty, S. and Hazarika, P.J. (2011). A survey of the theoretical developments in univariate skew normal distributions. *Assam Statistical Review*, 25(1), 41-63.

6. Chakraborty, S., Hazarika, P. J., and Ali, M. M. (2012). A New Skew Logistic Distribution and its Properties. *Pak. J. Statist*, 28(4), 513-524.

7. Chakraborty, S., Hazarika, P. J., and Ali, M. M. (2014). A Multimodal Skew Laplace Distribution. *Pak. J. Statist*, 30(2), 253-264.

8. Chakraborty, S., Hazarika, P. J., and Ali, M. M. (2015). A multimodal skewed extension of normal distribution: its properties and applications. *Statistics*, 49(4), 859-877.
9. Elal-Olivero, D. (2010). Alpha-skew-normal distribution. *Proyecciones (Antofagasta)*, 29(3), 224-240.
10. Harandi, S. S., and Alamatsaz, M. H. (2013). Alpha Skew Laplace distribution. *Statistics and Probability Letters*, 83(3), 774-782.
11. Hazarika, P. J., and Chakraborty, S. (2014). Alpha-Skew-Logistic Distribution. *IOSR Journal of Mathematics*, 10(4), 36-46.
12. Huang, W. J., and Chen, Y. H. (2007). Generalized skew-Cauchy distribution. *Statistics and Probability Letters*, 77(11), 1137-1147.
13. Nekoukhou, V., and Alamatsaz, M. H. (2012). A family of skew-symmetric-Laplace distributions. *Statistical Papers*, 53(3), 685-696.
14. Shafiei, S., Doostparast, M., and Jamalizadeh, A. (2016). The alpha–beta skew normal distribution: Properties and applications. *Statistics*, 50(2), 338-349.
15. Sharafi, M., Sajjadnia, Z., & Behboodian, J. (2017). A new generalization of alpha-skew-normal distribution. *Communications in Statistics-Theory and Methods*, 46(12), 6098-6111.
16. Venegas, O., Bolfarine, H., Gallardo, D. I., Vergara-Fernández, A., and Gómez, H. W. (2016). A Note on the Log-Alpha-Skew-Normal Model with Geochemical Applications. *Appl. Math*, 10(5), 1697-1703.
17. Wahed, A., and Ali, M. M. (2001). The skew-logistic distribution. *J. Statist. Res*, 35(2), 71-80.